\documentclass[12pt, a4paper]{amsart}
\usepackage{amscd, amssymb, pb-diagram}
\usepackage{amsthm}
\def\Aut{\operatorname{Aut}}

\def\range{\operatorname{range}}

\def\Ind{\operatorname{Ind}}
\def\Ad{\operatorname{Ad}}
\def\UHF{\operatorname{UHF}}
\def\C{\mathbb{C}}

\def\N{\mathbb{N}}
\def\Z{\mathbb{Z}}
\def\T{\mathbb{T}}
\def\TT{\mathcal{T}}
\def\LL{\mathcal{L}}
\def\OO{\mathcal{O}}
\def\KK{\mathcal{K}}

\def\HH{\mathcal{H}}
\newtheorem{thm}{Theorem}[section]
\newtheorem{cor}[thm]{Corollary}
\newtheorem{lemma}[thm]{Lemma}
\newtheorem{prop}[thm]{Proposition}

\theoremstyle{definition}
\newtheorem{definition}[thm]{Definition}
\theoremstyle{remark}
\newtheorem{remark}[thm]{Remark}
\newtheorem{example}[thm]{Example}

\numberwithin{equation}{section}
\begin{document}
\title[Exel's crossed Product and Relative Cuntz-Pimsner Algebras]{Exel's Crossed Product and\\Relative Cuntz-Pimsner Algebras}
\author{Nathan Brownlowe}

\author{Iain Raeburn}
\address{School of Mathematical and Physical Sciences\\
University of Newcastle\\
NSW 2308\\
Australia}
\email{nathan.brownlowe@studentmail.newcastle.edu.au}
\email{iain@maths.newcastle.edu.au}
\date{22 March 2004}
\begin{abstract}
We consider Exel's new construction of a crossed product of a $C^*$-algebra $A$ by an endomorphism $\alpha$. We prove that this crossed product is universal for an appropriate family of covariant representations, and we show that it can be realised as a relative Cuntz-Pimsner algbera. We describe a necessary and sufficient condition for the canonical map from $A$ into the crossed product to be injective, and present several examples to demonstrate the scope of this result. We also prove a gauge-invariant uniqueness theorem for the crossed product. 
\end{abstract}
\thanks{This research has been supported by the Australian Research Council}
\maketitle

\section{Introduction}
If $\alpha$ is an endomorphism of a $C^*$-algebra $A$, we can form a new $C^*$-algebra called the crossed product of $A$ by $\alpha$. This was first done by Cuntz  \cite{c}, and there are now several general theories \cite{pa,s,m}, which have been applied in a number of settings \cite{alnr,lr2,lr}.

In \cite{e1}, Exel proposed a new definition for the crossed product of a unital $C^*$-algebra $A$ by an endomorphism $\alpha$. Exel's crossed product depends not only on $A$ and $\alpha$, but also on the choice of a {\em transfer operator}, which is a positive continuous linear map $L:A\to A$ such that $L(\alpha(a)b)=aL(b)$ for $a,b\in A$. This new theory generalises previous constructions where the endomorphism is 
injective and has hereditary range \cite{m}, and has applications in the study
of classical irreversible dynamical systems \cite{ev}.

In this paper, we re-examine Exel's crossed product, denoted $A{\rtimes}_{\alpha,L}\N$, and identify a family of representations for which $A{\rtimes}_{\alpha,L}\N$ is universal. We then show that $A{\rtimes}_{\alpha,L}\N$ can be realised as a relative Cuntz-Pimsner algebra as in \cite{ms, fmr}, and use known results for relative Cuntz-Pimsner algebras to study $A{\rtimes}_{\alpha,L}\N$. In particular, we identify conditions which ensure that the canonical map $A\to A{\rtimes}_{\alpha,L}\N$ is injective, thus answering a question raised by Exel in \cite{e1}, and partially answered by him in \cite{e2}.

We begin with a brief discussion of relative Cuntz-Pimsner algebras, and we
state a lemma which we will use when considering the map 
$A\to A{\rtimes}_{\alpha,L}\N$. In \S\ref{Exel's Crossed Product} we discuss representations of Exel's crossed product. The main result in this section is the realization of $A{\rtimes}_{\alpha,L}\N$ as a relative Cuntz-Pimsner algebra. 

In \S\ref{The canonical map} we describe a necessary and sufficient condition on the transfer operator $L$ for $A\to A{\rtimes}_{\alpha,L}\N$ to be injective. We also show that this condition simplifies
when $A$ is a commutative $C^*$-algbera, and give examples to illustrate that our results do significantly improve those of Exel. In \S\ref{A Gauge-Invariant Uniqueness Theorem} we use our realisation of $A{\rtimes}_{\alpha,L}\N$ as a relative Cuntz-Pimsner algebra and results of Katsura \cite{k2} and Muhly-Tomforde \cite{mt} to prove a gauge-invariant uniqueness
theorem for $A{\rtimes}_{\alpha,L}\N$, which generalises the one of Exel and Vershik in \cite{ev}.
ÊÊÊÊÊ 
\section{Relative Cuntz-Pimsner algebras}\label{Relative Cuntz-Pimsner algebras}

Suppose that $A$ is a $C^*$-algebra and $X$ is a Hilbert bimodule over $A$, where the left action $a\cdot x$ is given by a homomorphism $\phi:A\to\LL(X)$, so that $a\cdot x=\phi(a)x$. A Toeplitz representation $(\psi,\pi)$ of $X$ in a $C^*$-algebra $B$ is a pair consisting of a linear map $\psi:X\to B$ and a homomorphism $\pi:A\to B$ such that
\[
\psi(x\cdot a)=\psi(x)\pi(a),\,{\psi(x)}^*\psi(y)=\pi({\langle x,y\rangle}_A),\text{ and }\psi(\phi(a)x)=\pi(a)\psi(x)
\]
for $x,y\in X$ and $a\in A$. Given such a representation, \cite[Proposition~1.6]{fr} says there is a homomorphism ${(\psi,\pi)}^{(1)}:\KK(X)\to B$ which satisfies
\[
{(\psi,\pi)}^{(1)}({\Theta}_{x,y})=\psi(x){\psi(y)}^*\text{ for } x,y\in X,
\]
and
\begin{equation}\label{a (1)-map condition}
{(\psi,\pi)}^{(1)}(T)\psi(x)=\psi(Tx)\text{ for }T\in\KK(X)\text{ and }x\in X.
\end{equation}
If $\rho:B\to C$ is a homomorphism of $C^*$-algebras, then $(\rho\circ\psi,\rho\circ\pi)$ is a Toeplitz representation of $X$, and we have
\[
{(\rho\circ\psi,\rho\circ\pi)}^{(1)}({\Theta}_{x,y})=\rho\circ\psi(x){\rho\circ\psi(y)}^*=\rho\circ{(\psi,\pi)}^{(1)}({\Theta}_{x,y})\text{ for all }x,y\in X.
\]
It follows from linearity and continuity that we have
\begin{equation}\label{pullrhothrough}
{(\rho\circ\psi,\rho\circ\pi)}^{(1)}= \rho\circ{(\psi,\pi)}^{(1)}.
\end{equation}

We define
\[
J(X):={\phi}^{-1}(\KK(X)),
\]
which is a closed two-sided ideal in $A$. Let $K$ be an ideal contained in $J(X)$. Following Muhly and Solel, we say that a Toeplitz representation $(\psi,\pi)$ of $X$ is {\em coisometric on} $K$ if
\begin{equation*}
{(\psi,\pi)}^{(1)}(\phi(a))=\pi(a)\quad\text{for all \ } a\in K.
\end{equation*}

\begin{prop}\label{RCP algebra}\cite[Proposition~1.3]{fmr}
Let $X$ be a Hilbert bimodule over $A$, and let $K$ be an ideal in $J(X)$. Then there are a $C^*$-algebra $\OO(K,X)$ and a Toeplitz representation $(k_X,k_A):X\to\OO(K,X)$ which is coisometric on $K$ and satisfies:
\begin{itemize}
\item[(i)] for every Toeplitz representation $(\psi,\pi)$ of $X$ which is coisometric on $K$, there is a homomorphism $\psi{\times}_K\pi$ of $\OO(K,X)$ such that $(\psi{\times}_K\pi)\circ k_X=\psi$ and $(\psi{\times}_K\pi)\circ k_A=\pi$; and
\item[(ii)] $\OO(K,X)$ is generated as a $C^*$-algebra by $k_X(X)\cup k_A(A)$.
\end{itemize}
The triple $(\OO(K,X),k_X,k_A)$ is unique: if $(B,k_X',k_A')$ has similiar properties, there is an isomorphism $\theta:\OO(K,X)\to B$ such that $\theta\circ k_X=k_X'$ and $\theta\circ k_A=k_A'$. There is a strongly continuous gauge action $\gamma:\T\to\Aut\OO(K,X)$ which satisfies ${\gamma}_z(k_A(a))=k_A(a)$ and ${\gamma}_z(k_X(x))=zk_X(x)$ for $a\in A,x\in X$.
\end{prop}
Ê 
The algebra $\OO(K,X)$ is called the relative Cuntz-Pimsner algebra determined by~$K$, and was first studied by Muhly and Solel in \cite{ms}. The algebra $\OO(\{0\},X)$ is the Toeplitz algebra $\TT(X)$ (see \cite[Proposition~1.4]{fr}), and $\OO(J(X),X)$ is the Cuntz-Pimsner algebra $\OO(X)$ \cite{p}.
The following lemma tells us when $k_A:A\to\OO(K,X)$ is injective. 

\begin{lemma}\label{injectivity of k_A}
Let $X$ be a Hilbert bimodule over $A$ and let $(\OO(K,X),k_A,k_X)$ be a relative Cuntz-Pimsner algebra associated to $X$. Then $k_A$ is injective if and only if ${\phi}|_K:K\to\LL(X)$ is injective.
\end{lemma}
\begin{proof}
If ${\phi}|_K$ is injective, then \cite[Proposition~2.21]{ms} implies that $k_A$ is injective. Conversely, suppose $k_A$ is injective and $a\in K$ satisfies ${\phi}|_K(a)=0$. Then $k_A(a)={(k_X,k_A)}^{(1)}({\phi}|_K(a))=0$, and since $k_A$ is injective, this implies $a=0$. Thus ${\phi}{|}_K:K\to\LL(X)$ is injective.
\end{proof}

\section{Exel's Crossed Product}\label{Exel's Crossed Product}
Let $A$ be a unital $C^*$-algebra and $\alpha$ an endomorphism of $A$; we do not assume that $\alpha$ is unital or injective. In \cite{e1}, Exel defined a \emph{transfer operator} $L$ for $(A,\alpha)$ to be a continuous linear map $L:A\to A$ such that
\begin{itemize}
\item[(i)] $L$ is positive in the sense that $a\ge 0\Longrightarrow L(a)\ge0$, and
\item[(ii)] $L(\alpha(a)b)=aL(b)$, for all $a,b\in A$.
\end{itemize}
He then defined $\TT(A,\alpha,L)$ to be the universal unital $C^*$-algebra generated by a copy of $A$ and an element $S$ satisfying the relations $Sa=\alpha(a)S$ and $S^*aS=L(a)$ for $a\in A$, so that $\TT(A,\alpha,L)$ is by definition universal for the following representations. 
\begin{definition}
A pair $(\rho,V)$, consisting of a unital homomorphism $\rho$ of $A$ into a $C^*$-algebra $B$ and an element $T\in B$, is a \emph{Toeplitz-covariant representation of} $(A,\alpha,L)$ in $B$ if for every $a\in A$,
\begin{itemize}
\item[(TC1)] $V\rho(a)=\rho(\alpha(a))V$, and 
\smallskip
\item[(TC2)] $V^*\rho(a)V=\rho(L(a))$. 
\end{itemize}
We denote by $(i_A,S)$, the universal Toeplitz-covariant representation of $(A,\alpha,L)$ in $\TT(A,\alpha,L)$. If $(\rho, V)$ is a Toeplitz-covariant representation of $(A,\alpha,L)$, we denote by $\rho\times V$ the representation of $\TT(A,\alpha,L)$ such that $(\rho\times V)\circ i_A=\rho$ and $(\rho\times V)(S)=V$.
\end{definition}
 
The homomorphism $i_A:A\to \TT(A,\alpha,L)$ is injective: to see this, we need an example of a Toeplitz-covariant representation $(\rho,V)$ with $\rho$ injective, and one such example is given in \cite{e1}.

Given the triple $(A,\alpha,L)$, we recall from \cite{e1} the construction of the Hilbert $A$-bimodule $M_L$. We let $A_L$ be a copy of the underlying vector space of $A$. We define a right action of $A$ on $A_L$ by 
\[
m\cdot a=m\alpha(a)\ \text{for $m\in A_L$ and $a\in A$,}
\] 
and an $A$-valued map ${\langle\cdot,\cdot\rangle}_L$ on $A_L$ by
\[
{\langle m,n\rangle}_L=L(m^*n)\ \text{for $m,n\in A_L$.}
\]
We define $N:=\{a\in A_L:{\langle a,a\rangle}_L=0\}$; it follows from the Cauchy-Schwarz inequality that $N$ is a subspace of $A_L$, and we can form the quotient space $A_L/N$. We denote the quotient map by $q:A_L\to A_L/N$, and then $A_L/N$ is a right $A$-module with inner-product ${\langle q(a),q(b)\rangle}_L=L(a^*b)$. By completing $A_L/N$ we get a right Hilbert $A$-module which we denote by $M_L$. 
For $a\in A$ and $m\in A_L$ we have
\[
\|\langle am,am\rangle_L\|=\|L(m^*a^*am)\|\le{\|a\|}^2\|L(m^*m)\|={\|a\|}^2\|\langle m,m\rangle_L\|,
\]
and it follows that left multiplication by $a$ on $A_L$ extends to a bounded adjointable operator on $M_L$. ThisÊ defines a homomorphism $\phi:A\to\LL(M_L)$, and writing $\phi(a)m:=a\cdot m$ makes $M_L$ a Hilbert bimodule over $A$. Note that $q(A_L)$ is dense in $M_L$.

In the following lemma we see that there is a one-to-one correspondence between Toeplitz-covariant representations of $(A,\alpha,L)$ and Toeplitz representations of $M_L$.

\begin{lemma}\label{induced T-R}
Given a Toeplitz-covariant representation $(\rho,V)$ of $(A,\alpha,L)$ in a $C^*$-algebra $B$, there exists a linear map ${\psi}_V:M_L\to B$ such that ${\psi}_V(q(a))=\rho(a)V$ and the pair $({\psi}_V,\rho)$ is a Toeplitz representation of $M_L$ in $B$. Conversely, if $(\psi,\pi)$ is a Toeplitz representation of $M_L$ in $B$ and $\pi$ is unital, then the pair $(\pi,\psi(q(1)))$ is a Toeplitz-covariant representation of $(A,\alpha,L)$, and $\psi_{\psi(q(1))}=\psi$.
\end{lemma}
\begin{proof}
We define ${\theta}:A_L\to B$ by ${\theta}(a)=\rho(a)V$. Then $\theta$ is linear, and for $a\in A$ we have
\begin{align*}
{\|{\theta}(a)\|}^2 &= {\|\rho(a)V\|}^2= \|{(\rho(a)V)}^*\rho(a)V\|= \|V^*\rho(a^*a)V\|= \|\rho(L(a^*a))\|\\
&\le \|L(a^*a)\|= \|{\langle a,a\rangle}_L\|,
\end{align*}
so $\theta$ is bounded for the semi-norm on $A_L$. Thus ${\theta}$ induces a bounded map ${\psi}_V:M_L\to B$ satisfying ${\psi}_V(q(a))=\rho(a)V$ for $a\in A$.
For $a,b,c\in A$ we have
\begin{gather*}
{\psi}_V(q(b)\cdot a)={\psi}_V(q(b)\alpha(a))=\rho(b\alpha(a))V=\rho(b)V\rho(a)={\psi}_V(q(b))\rho(a),\\
{{\psi}_V(q(b))}^*{\psi}_V(q(c))={(\rho(b)V)}^*\rho(c)V=V^*\rho(b^*c)V=\rho(L(b^*c))=\rho({\langle q(b),q(c)\rangle}_L),\mbox{ and }\\
{\psi}_V(a\cdot q(b))={\psi}_V(aq(b))=\rho(ab)V=\rho(a)\rho(b)V=\rho(a){\psi}_V(q(b)).
\end{gather*}
Thus $({\psi}_V,\rho)$ is a Toeplitz representation of $M_L$ in $B$.

Now let $(\psi,\pi):M_L\toÊ B$ be a Toeplitz representation of $M_L$ inÊa $C^*$-algebra $B$ with $\pi$ unital. Then for $a\in A$ we have
\[
\psi(q(1))\pi(a)=\psi(q(1)\cdot a)=\psi(q(\alpha(a))=\psi(\alpha(a)\cdot q(1))=\pi(\alpha(a))\psi(q(1)),
\]
and
\begin{align*}
{\psi(q(1))}^*\pi(a)\psi(q(1)) &= {\psi(q(1))}^*\psi(a\cdot q(1))={\psi(q(1))}^*\psi(q(a))\\
&= \pi({\langle q(1),q(a)\rangle}_L)= \pi(L(1^*a))=\pi(L(a)),
\end{align*}
so $(\pi,\psi(q(1)))$ is a Toeplitz-covariant representation of $(A,\alpha,L)$. Finally, for $a\in A$ we have
\[
\psi_{\psi(q(1))}(q(a))=\pi(a)\psi(q(1))=\psi(a\cdot q(1))=\psi(q(a)),
\]
which implies that $\psi_{\psi(q(1))}=\psi$.
\end{proof}

\begin{cor}\label{the Toeplitz algebra}
The $C^*$-algbera $\TT(A,\alpha,L)$ is isomorphic to the Toeplitz algebra $\TT(M_L)$.
\end{cor}

\begin{proof}
We prove that $\TT(A,\alpha, L)$ has the universal property which characterises $\TT(M_L)$. Applying the lemma to the pair $(i_A,S)$ gives a Toeplitz representation $(\psi_S,i_A)$ of $M_L$ in $\TT(A,\alpha, L)$, which generates $\TT(A,\alpha, L)$ because $i_A$ and $S$ do. Now suppose $(\psi,\pi)$ is a Toeplitz representation of $M_L$. Note that $M_L$ is essential as a left $A$-module, in the sense that $A\cdot M_L=M_L$. This implies that the essential subspace $\pi(1)\HH$ is reducing for $(\psi,\pi)$, so we can apply the lemma to the restriction of $(\psi,\pi)$ to $\pi(1)\HH$; this gives a Toeplitz-covariant representation $(\pi|,\psi|(q(1)))$ on $\pi(1)\HH$. Now the representation $\mu:=(\pi|\times \psi|(q(1)))\oplus 0$ has $\mu\circ i_A=\pi|\oplus 0=\pi$, and for $a\in A$ we have
\begin{align*}
\mu\circ\psi_S(q(a))&=\mu(i_A(a)S)=\mu(i_A(a))\mu(S)=(\pi|(a)\psi|(q(1)))\oplus 0\\
&=\pi(a)\psi(q(1))=\psi_{\psi(q(1))}(q(a))=\psi(q(a)),
\end{align*}
which implies that $\mu\circ\psi_S=\psi$. 
\end{proof}

Corollary~\ref{the Toeplitz algebra} has been obtained independently by Nadia Larsen. 

\begin{remark}\label{(1)-map}
The Toeplitz representation $({\psi}_S,i_A)$ induces a homomorphism ${({\psi}_S,i_A)}^{(1)}$ of $\KK(M_L)$ into $\TT(A,\alpha,L)$. We claim that ${({\psi}_S,i_A)}^{(1)}$ is injective. To see this, let $\pi:\TT(A,\alpha,L)\to B(\HH)$ be a faithful non-degenerate representation of $\TT(A,\alpha,L)$. Then, as in the proof of \cite[Proposition~1.6]{fr}, we have ${({\psi}_S,i_A)}^{(1)}:={\pi}^{-1}\circ\Ad U\circ\Ind(\pi\circ i_A)$, where $U:M_L{\otimes}_A\HH\to\HH$ is an isometry given by $U(m{\otimes}_Ah)=\pi({\psi}_S(m))h$. Since $\pi\circ i_A$ is faithful, the induced representation is faithful \cite[Corollary~2.74]{rw}, and ${({\psi}_S,i_A)}^{(1)}$ is injective, as claimed. 

The range of any homomorphism of $C^*$-algebras is closed, and since ${({\psi}_S,i_A)}^{(1)}(\KK(M_L))$ is dense in $\overline{{\psi}_S(M_L){{\psi}_S(M_L)}^*}$, it follows that ${({\psi}_S,i_A)}^{(1)}$ is an isomorphism of $\KK(M_L)$ onto the $C^*$-algebra $\overline{{\psi}_S(M_L){{\psi}_S(M_L)}^*}=\overline{i_A(A)SS^*i_A(A)}$.
\end{remark}

We will now discuss Exel's notion of a redundancy. Define $M:=\overline{i_A(A)S}={\psi}_S(M_L)$. Conditions (TC1) and (TC2) imply that $i_A(A)M\subseteq M$, $Mi_A(A)\subseteq M$ and $M^*M\subseteq i_A(A)$, so $M$ is a Hilbert bimodule over $i_A(A)$. It follows that left multiplication by elements of $i_A(A)$ on $M$ could coincide with left multiplication by elements in $\overline{MM^*}=\overline{i_A(A)SS^*i_A(A)}$. In \cite{e1}, Exel defines a {\em redundancy} to be a pair $(i_A(a),k)$ such that $a\in A$, $k\in\overline{i_A(A)SS^*i_A(A)}$ and 
\[
i_A(a)i_A(b)S=ki_A(b)S\,\text{ for all }b\in A.
\]
The next lemma provides a useful identification of the redundancies.

\begin{lemma}\label{our redundancies}
Let $a\in A$ and let $k\in \TT(A,\alpha, L)$. Then $(i_A(a),k)$ is a redundancy if and only if $a\in J(M_L):=\phi^{-1}(\KK(M_L))$ and $k={({\psi}_S,i_A)}^{(1)}(\phi(a))$.
\end{lemma}

\begin{proof}
First suppose that $a\in J(M_L)$ and $k={({\psi}_S,i_A)}^{(1)}(\phi(a))$. Then $k$ belongs to the image $\overline{i_A(A)SS^*i_A(A)}$ of ${({\psi}_S,i_A)}^{(1)}$, and for $b\in A$ we have 
\begin{align*}
i_A(a)i_A(b)S &= i_A(a){\psi}_S(q(b))= {\psi}_S(\phi(a)q(b))\\
&={({\psi}_S,i_A)}^{(1)}(\phi(a)){\psi}_S(q(b))\\
&= {({\psi}_S,i_A)}^{(1)}(\phi(a))i_A(b)S,
\end{align*}
where the second last equality follows from Equation~(\ref{a (1)-map condition}). Thus $(i_A(a),k)$ is a redundancy.

Now suppose that $(i_A(a),k)$ is a redundancy. It follows from Remark~\ref{(1)-map} that there exists a unique $t\in\KK(M_L)$ such that ${({\psi}_S,i_A)}^{(1)}(t)=k$. Then for $b\in A$ we have
\begin{align*}
{\psi}_S\big(\phi(a)(q(b))\big) &= {\psi}_S(q(ab))= i_A(ab)S= i_A(a)i_A(b)S\\
&= ki_A(b)S= {({\psi}_S,i_A)}^{(1)}(t){\psi}_S(q(b))= {\psi}_S(t(q(b))),
\end{align*}
Since $i_A:A\to \TT(A,\alpha,L)$ is injective, ${\psi}_S$ is also injective, and it follows that $\phi(a)(m)=t(m)$ for all $m\in M_L$. Hence $\phi(a)=t$, and the result follows. 
\end{proof}

Exel defined the crossed product of $(A,\alpha,L)$ to be the quotient of $\TT(A,\alpha,L)$ by the ideal generated by the set 
\[
\{i_A(a)-k:(i_A(a),k)\text{ is a redundancy with }a\in\overline{A\alpha(A)A}\}.
\]
We denote the quotient map by $Q:\TT(A,\alpha,L)\to A{\rtimes}_{\alpha,L}\N$. The next corollary follows immediately from Lemma~\ref{our redundancies}.

\begin{cor}\label{the crossed product}
Let $K_{\alpha}:=\overline{A\alpha(A)A}\cap J(M_L)$
and denote by $I(A,\alpha,L)$ the ideal in $\TT(A,\alpha,L)$ generated by 
\[
\{i_A(a)-{({\psi}_S,i_A)}^{(1)}(\phi(a)):a\in K_{\alpha}\}.
\]
Then $A{\rtimes}_{\alpha,L}\N$ is $\TT(A,\alpha,L)/I(A,\alpha,L)$. 
\end{cor}

To describe $A{\rtimes}_{\alpha,L}\N$ as a universal object, we need to identify the Toeplitz-covariant representations that vanish on the ideal $I(A,\alpha,L)$. We need a lemma:

\begin{lemma}\label{equ for redundancies}
Suppose $(\rho,V)$ is a covariant representation of $(A,\alpha,L)$. Then we have
\begin{equation}\label{vanishing repns}
(\rho\times V)\circ{({\psi}_S,i_A)}^{(1)}={({\psi}_V,\rho)}^{(1)}.
\end{equation}
\end{lemma}

\begin{proof}
We know from \eqref{pullrhothrough} that
\[
(\rho\times V)\circ{({\psi}_S,i_A)}^{(1)}={((\rho\times V)\circ{\psi}_S,(\rho\times V)\circ i_A)}^{(1)}.
\]
Since $(i_A,S)$ is the universal Toeplitz-covariant representation, we have $(\rho\times V)\circ i_A=\rho$, and $(\rho\times V)(S)=V$. So for $a\in A$ we have
\[
(\rho\times V)\circ{\psi}_S(q(a))=\rho\times V(i_A(a)S)=\rho(a)V=\psi_V(q(a)),
\]
 and hence we also have $(\rho\times V)\circ{\psi}_S=\psi_V$.
\end{proof}

Equation~(\ref{vanishing repns}) motivates the following definition. 

\begin{definition}\label{cov repn}
Consider the triple $(A,\alpha,L)$, and let $(\rho,V)$ be a Toeplitz-covariant representation in a $C^*$-algebra $B$. We say that $(\rho,V)$ is a \emph{covariant representation} of $(A,\alpha,L)$ if in addition we have 
\begin{itemize}
\item[(C3)] $\rho(a)={({\psi}_V,\rho)}^{(1)}(\phi(a))$ for all $a\in K_{\alpha}$.
\end{itemize}
\end{definition}
The following Proposition says that $A{\rtimes}_{\alpha,L}\N$ is universal for covariant representations of $(A,\alpha,L)$.

\begin{prop}\label{uni prop}
Let $\alpha$ be an endomorphism of a unital $C^*$-algbera $A$, and let $L$ be a transfer operator for $(A,\alpha)$. The pair $(j_A,T):=(Q\circ i_A,Q(S))$ is a covariant representation of $(A,\alpha,L)$ in $A{\rtimes}_{\alpha,L}\N$, and for every covariant representation $(\rho,V)$ of $(A,\alpha,L)$, there is a representation ${\tau}_{\rho,V}$ of $A{\rtimes}_{\alpha,L}\N$ such that ${\tau}_{\rho,V}\circ j_A=\rho$ and ${\tau}_{\rho,V}(T)=V$.
\end{prop}

\begin{proof}
The pair $(Q\circ i_A,Q(S))$ is Toeplitz-covariant because $(i_A,S)$ is, and its integrated form $(Q\circ i_A)\times Q(S)$ is precisely $Q$. By Lemma~\ref{induced T-R}, we get a Toeplitz representation $({\psi}_{Q(S)},Q\circ i_A):M_L\to A{\rtimes}_{\alpha,L}\N$, and for $a\in K_{\alpha}$ we have
\begin{align*}
(Q\circ i_A)(a) &= Q(i_A(a))=Q\big({({\psi}_S,i_A)}^{(1)}(\phi(a))\big)\\
&= ((Q\circ i_A)\times Q(S))\big({({\psi}_S,i_A)}^{(1)}(\phi(a))\big)\\
&= {({\psi}_{Q(S)},Q\circ i_A)}^{(1)}(\phi(a)),
\end{align*}
using Lemma~\ref{equ for redundancies}. So the pair $(Q\circ i_A,Q(S))$ is covariant.

Now suppose $(\rho,V)$ is a covariant representation of $(A,\alpha,L)$. The Toeplitz-covariant representation $(\rho,V)$ gives us a representation $\rho\times V$ of $\TT(A,\alpha,L)$, and condition (C3) says that $\rho\times V$ vanishes on the generators of the ideal $I(A,\alpha,L)$. Hence Corollary~\ref{the crossed product} implies that $\rho\times V$ factors through a representationÊ ${\tau}_{\rho,V}$ of $A{\times}_{\alpha,L}\N$. Then
\begin{gather*}
{\tau}_{\rho,V}\circ j_A={\tau}_{\rho,V}\circ Q\circ i_A=(\rho\times V)\circ i_A=\rho,\text{ and}\\
{\tau}_{\rho,V}(T)={\tau}_{\rho,V}(Q(T))=(\rho\times V)(T)=V,
\end{gather*}
so ${\tau}_{\rho,V}$ has the required properties.
\end{proof}

We now realise $A{\rtimes}_{\alpha,L}\N$ as a relative Cuntz-Pimsner algebra. 
 
\begin{prop}\label{isom of Cross P and rel CP algebra}
Suppose $\alpha$ is an endomorphism of a unital $C^*$-algebra $A$ and $L$ is a transfer operator for $(A,\alpha)$. Then there is an isomorphism $\theta: \OO(K_\alpha,M_L)\to A{\rtimes}_{\alpha,L}\N$ such that $\theta\circ k_A=j_A$ and $\theta(k_{M_L}(q(1)))=T$. 
\end{prop}
 
\begin{proof}
Consider the triple $(A,{\psi}_T,j_A)$, where $({\psi}_T,j_A)$ is the Toeplitz representation of $M_L$ induced by the pair $(j_A,T)$, as in Lemma~\ref{induced T-R}.  
We will prove that $(A{\rtimes}_{\alpha,L}\N,{\psi}_T,j_A)$ satisfies the conditions of Proposition~\ref{RCP algebra}.

Since $(j_A,T)$ is covariant, it satisfies (C3), which says precisely that  $(\psi_T,j_A)$ is coisometric on $K_\alpha$. Let $(\psi,\pi)$ be a Toeplitz representation of $M_L$ which is coisometric on $K_{\alpha}$; since $M_L$ is essential, we suppose by throwing away a trivial representation that $\pi$ is unital (see the proof of Corollary~\ref{the Toeplitz algebra}). Then Lemma~\ref{induced T-R} gives a Toeplitz-covariant representation $(\pi,\psi(q(1)))$. Since ${\psi}_{\psi(q(1))}=\psi$ and $(\psi,\pi)$ is coisometric on $K_\alpha$, $(\pi,\psi(q(1)))$ is covariant. Now Proposition~\ref{uni prop} gives a representation $\tau_{\pi,\psi(q(1))}$ of $A{\rtimes}_{\alpha,L}\N$ such that $\tau_{\pi,\psi(q(1))}\circ j_A=\pi$ and $\tau_{\pi,\psi(q(1))}(T)=\psi(q(1))$. For $a\in A$ we have
\[
\tau_{\pi,\psi(q(1))}(\psi_T(q(a)))=\tau_{\pi,\psi(q(1))}(j_A(a)T)=\pi(a)\psi(q(1))=\psi(q(a)),
\]
and it follows that $\tau_{\pi,\psi(q(1))}\circ \psi_T=\psi$. So $\psi\times_{K_\alpha}\pi:=\tau_{\pi,\psi(q(1))}$ satisfies condition (i) of Proposition~\ref{RCP algebra}. Since $\psi_T(M_L)\cup j_A(A)$ generates $A{\rtimes}_{\alpha,L}\N$, condition (ii) is also satisfied, and applying Proposition~\ref{RCP algebra} gives the result.
\end{proof}

Notice that when $\alpha(1)=1$, we have $K_{\alpha}=J(M_L)$, and the crossed product $A{\rtimes}_{\alpha,L}\N$ is the Cuntz-Pimsner algebra $\OO(M_L)$.
\section{Injectivity of $j_A:A\to A{\rtimes}_{\alpha,L}\N$}\label{The canonical map}

\begin{definition}\label{almost faithful}
Suppose that $A$ is a unital $C^*$-algebra, $\alpha$ is an endomorphism of $A$ and $L$ is a transfer operator for $(A,\alpha)$. We say that $L$ is \emph{faithful} on an ideal $I$ of $A$ if
\[
a\in I\mbox{ and }L(a^*a)=0\Longrightarrow a=0;
\]
we say that $L$ is \emph{almost faithful} on $I$ if
\[
a\in I\mbox{ and }L({(ab)}^*ab)=0 \text{ for all } b\in A\Longrightarrow a=0.
\]
\end{definition}

\begin{thm}\label{injectivity of j_A for non-comm}
Let $\alpha$ be an endomorphism of a unital $C^*$-algebra $A$, and let $L$ be a transfer operator for $(A,\alpha)$. Then the map $j_A:A\to A{\rtimes}_{\alpha,L}\N$ is injective if and only if $L$ is almost faithful on $K_{\alpha}=\overline{A\alpha(A)A}\cap J(M_L)$.
\end{thm}

\begin{proof}
It follows from Proposition~\ref{isom of Cross P and rel CP algebra} that the map $j_A$ is injective if and only if $k_A:A\to\OO(K_{\alpha},M_L)$ is injective. By Lemma~\ref{injectivity of k_A} this is true if and only if ${\phi}|_{K_{\alpha}}:K_{\alpha}\to \LL(M_L)$ is injective, and so it suffices to prove that the transfer operator $L$ is almost faithful on $K_{\alpha}$ if and only if ${\phi}|_{K_{\alpha}}:K_{\alpha}\to\LL(M_L)$ is injective. But for $a\in K_\alpha$ and $b\in A$, we have
\begin{align*}
\|L({(ab)}^*ab)\| &= \|{\langle q(ab),q(ab)\rangle}_L\|= {\|q(ab)\|}^2= {\|a\cdot q(b)\|}^2\\
&= {\|\phi(a)(q(b))\|}^2= {\|{\phi}|_{K_{\alpha}}(a)(q(b))\|}^2,
\end{align*}
and this implies the desired equivalence. 
\end{proof}

\begin{cor}\label{injectivity of j_A}
Let $\alpha$ be an endomorphism of a unital commutative $C^*$-algebra $A$, and let $L$ be a transfer operator for $(A,\alpha)$. Then the map $j_A:A\to A{\rtimes}_{\alpha,L}\N$ is injective if and only if $L$ is faithful on $K_{\alpha}$.
\end{cor}
\begin{proof}
If $L$ is faithful on $K_{\alpha}$ then it follows from Theorem~\ref{injectivity of j_A for non-comm} that $j_A:A\to A{\rtimes}_{\alpha,L}\N$ is injective. Conversely, suppose $j_A:A\to A{\rtimes}_{\alpha,L}\N$ is injective. By Theorem~\ref{injectivity of j_A for non-comm}, this implies that $L$ is almost faithful on $K_{\alpha}$. Suppose $a\in K_{\alpha}$ satisfies $L(a^*a)=0$. Then for every $b\in A$ we have
\[
\|L({(ab)}^*ab)\|=\|L({(ba)}^*ba)\|=\|L(a^*b^*ba)\|\le {\|b\|}^2\|L(a^*a)\|=0.
\]
Thus $L({(ab)}^*ab)=0$ for every $b\in A$, which implies $a=0$, and we have shown that $L$ is faithful on $K_{\alpha}$.
\end{proof}

In \cite{e2}, Exel assumed that $\alpha$ is a unital injective endomorphism and $L={\alpha}^{-1}\circ E$, where $E$ is a conditional expectation of $A$ onto $\alpha(A)$ satisfying $E(a^*a)=0\Longrightarrow a=0$ (Exel says $E$ is non-degenerate). Under these conditions he proves that $j_A:A\to A{\rtimes}_{\alpha,L}\N$ is injective \cite[Theorem~4.12]{e2}. Notice that such $L$ are faithful, and so \cite[Theorem~4.12]{e2} follows from Theorem~\ref{injectivity of j_A for non-comm}. The following examples  show that our theorem is stronger in several different ways.

\begin{example}\label{the backward and forward shift} In this example, the endomorphism is not unital.
Let $A=c$, the space of convergent sequences under the $\sup$ norm, and let $\alpha$ be the forward shift ${\tau}_f$.  Then the backward shift $L={\tau}_b$ is a transfer operator for $(c,{\tau}_f)$ and we have
\[
M_{{\tau}_b} = c/\C e_0,\,J(M_{{\tau}_b}) = c,\text{ and }K_{{\tau}_f} = \{f\in c:f(0)=0\};
\]
notice that $L={\tau}_b$ is faithful on $K_{{\tau}_f}$, but not on all of $c$. It follows from Corollary~\ref{injectivity of j_A} that the map $j_c:c\to c{\rtimes}_{{\tau}_f,{\tau}_b}\N$ is injective.
\end{example}

\begin{example}
In this example, the endomorphism is not injective. Again the algebra $A$ is $c$, but now
we view the backward shift ${\tau}_b$ as the endomorphism, and take for $L$ the forward shift ${\tau}_f$. Then we have $M_{{\tau}_f}=A_L$, $J(M_{{\tau}_f}) = c$, and $K_{{\tau}_b} = c$. In this case, $L={\tau}_f$ is faithful on $K_{{\tau}_b}$, so Corollary~\ref{injectivity of j_A} shows that $j_c:c\to c{\rtimes}_{{\tau}_b,{\tau}_f}\N$ is injective.
\end{example}

\begin{example}\label{uhf} In this example, the transfer operator is almost faithful but is not faithful. We take $A$ to be the UHF algebra $\UHF(n^\infty)$, viewed as the direct limit $\varinjlim(A_N,i_N)$ with $A_N=\bigotimes_{k=1}^NM_n(\C)$ and 
\[
i_N(a_1\otimes\dots\otimes a_N):=a_1\otimes\dots\otimes a_N\otimes 1;
\]
we denote the canonical embeddings by $i^N:A_N\to A$. The mapsÊ ${\alpha}_N:A_N\to A_{N+1}$ defined by 
\[
{\alpha}_N(a_1\otimes\dots\otimes a_N)=e_{11}\otimes a_1\otimes\dots\otimes a_N,
\]
induce an injective endomorphism $\alpha:A\to A$ such that $\alpha(i^N(a))=i^{N+1}({\alpha}_N(a))$ for $a\in A_N$. Since $\range\alpha$ is closed, it follows that $\range\alpha=i^1(e_{11})A\,i^1(e_{11})$. We can then define $L:A\to A$ by
\[
L(a)={\alpha}^{-1}(i^1(e_{11})a\,i^1(e_{11})).
\]
Then $L$ is positive, continuous and linear. To see that $L$ is a transfer operator, let $a=\bigotimes a_i\in A_N$, $b=\bigotimes b_i\in A_{N+1}$, and compute:
\begin{align*}
L(\alpha(i^N(a))i^{N+1}(b)) &= L(i^{N+1}(e_{11}b_1\otimes a_1b_2\otimes\dots\otimes a_Nb_{N+1}))\\
&= {\alpha}^{-1}(i^1(e_{11})i^{N+1}(e_{11}b_1\otimes a_1b_2\otimes\dots\otimes a_Nb_{N+1})i^1(e_{11}))\\
&= {(b_1)}_{11}{\alpha}^{-1}(i^{N+1}(e_{11}\otimes a_1b_2\otimes\dots\otimes a_Nb_{N+1}))\\
&= {(b_1)}_{11}i^{N}(a_1b_2\otimes\dots\otimes a_Nb_{N+1})\\
&= {(b_1)}_{11}i^N(a_1\otimes\dots\otimes a_N)i^{N}(b_2\otimes\dots\otimes b_{N+1})\\
&= i^N(a){(b_1)}_{11}{\alpha}^{-1}(i^{N+1}(e_{11}\otimes b_2\otimes\dots\otimes b_M))\\
&= i^N(a){\alpha}^{-1}(i^1(e_{11})i^{N+1}(b)i^1(e_{11}))\\
&= i^N(a)L(i^{N+1}(b)).
\end{align*}
It follows from linearity and continuity of $L$ and $\alpha$ that $L(\alpha(a)b)=aL(b)$ for all $a,b\in A$, and hence $L$ is a transfer operator for $(A,\alpha)$.

For $j\in\{1,\dots,n\}$ define $b_j:=i^1(e_{j1})$. Suppose $a\in A$ satisfies $L({(ab)}^*ab)=0$ for all $b\in A$. Then $0=L({(ab_j)}^*ab_j)={\alpha}^{-1}(i^1(e_{11}){b_j}^*a^*ab_ji^1(e_{11}))$ for all $j$, and this implies that $ab_j=0$ for all $j$. Thus
\[
0=\sum_{j=1}^nab_j{b_j}^*=ai^1\Big(\sum_{j=1}^ne_{jj}\Big)=ai^1(1)=a,
\]
and hence $L$ is almost faithful on $A$.
To see that $L$ is not faithful we let $a_0\in M_n(\C)$ be a non-zero matrix whose first column is zero. Then ${({a_0}^*a_0)}_{11}=0$ and 
\begin{align*}
L({i^1(a_0)}^*i^1(a_0)) = {\alpha}^{-1}(i^1(e_{11}{a_0}^*a_0e_{11}))= {\alpha}^{-1}({({a_0}^*a_0)}_{11}i^1(e_{11}))= {\alpha}^{-1}(0)= 0,
\end{align*}
whereas $i^1(a_0)\not=0$ because $i^1$ is injective. 

The endomorphism $\alpha$ is injective and has hereditary range. Under these assumptions, Exel proved in \cite[Theorem~4.7]{e1} that $A{\rtimes}_{\alpha,L}\N$ is isomorphic to the Stacey crossed product $A{\rtimes}_{\alpha}\N$. This crossed product was first considered by Cuntz, who showed in \cite{c} that $\UHF(n^{\infty}){\rtimes}_{\alpha}\N$ is isomorphic to the Cuntz algebra $\OO_n$.
\end{example}

\begin{example} This is an example of a commutative $C^*$-algebra with a transfer operator $L$ which is not faithful on $K_{\alpha}$, so that $A$ does not embed in Exel's crossed product. Let $A:=C([0,2])$, and define $\alpha:C([0,2])\to C([0,2])$ by 
\[
\alpha(f)(x):=
\begin{cases}
f(2x)Ê 
ÊÊÊ & \text{if $x\in[0,1]$} \\
f(4-2x)
ÊÊÊ & \text{if $x\in(1,2]$.} \\
\end{cases}
\]
Then the map $L:C([0,2])\to C([0,2])$ defined by $L(f)(x)=f(x/2)$, is a transfer operator for $(A,\alpha)$.
We have $A_L=C([0,2])$ as a vector space, and 
\begin{align*}
N &:= \{f\in C([0,2]):L(f^*f)=0\}\\
&= \{f\in C([0,2]): f(x)=0 \text{ for all } x\in[0,1]\}.
\end{align*}
Thus the restriction map $r:f\mapsto f|_{0,1]}$ induces a vector-space isomorphism of $A_L/N$ onto $C([0,1])$, which converts the bimodule structure into
\[
\langle g,h\rangle_L(x)=\overline{g(x/2)}h(x/2),\ g\cdot f(x)=g(x)f(2x), \ f\cdot g(x)=f(x)g(x)
\]
for $g,h\in C([0,2])$ and $f\in A=C([0,2])$; it follows from the first formula that $r$ is isometric for the sup-norm on $C([0,1])$, so $A_L/N$ is complete and $M_L=A_L/N$. Now for $f\in A$ and $x\in [0,1]$, we have
\[
\Theta_{r(f),1}(g)(x)=r(f)(x)\langle 1,g\rangle_L(2x)=f(x)g(x)=(\phi(f)g)(x),
\]
so $f\in J(M_L)$. Thus $J(M_L)=A$, which implies $K_\alpha=A$ because $\alpha(1)=1$. The transfer function $L$ is not faithful on $C([0,2])$: any nonzero function $f\in C([0,2])$ with $f|_{[0,1]}=0$ will satisfy $L(f^*f)=0$. Hence it follows from Corollary~\ref{injectivity of j_A} that the canonical map $C([0,2])\to C([0,2]){\rtimes}_{\alpha,L}\N$ is not injective.
\end{example}

\section{Gauge Invariant Uniqueness Theorem}\label{A Gauge-Invariant Uniqueness Theorem}

Using the isomorphism $\theta:\OO(K_\alpha,M_L)\to A{\rtimes}_{\alpha,L}\N$ of Proposition~\ref{isom of Cross P and rel CP algebra}, we can see that there is a natural gauge action $\delta:\T\to \Aut (A{\rtimes}_{\alpha,L}\N)$ such that ${\delta}_z(j_A(a))=j_A(a)$, ${\delta}_z(T)=zT$ and $\theta\circ \gamma_z=\delta_z\circ \theta$. 

\begin{thm}\label{the giut}
Let $\alpha$ be an endomorphism of a unital $C^*$-algebra $A$, and let $L$ be a transfer operator for $(A,\alpha)$. Suppose $B$ is a $C^*$-algebra and $(\rho,V)$ is a covariant representation of $(A,\alpha,L)$ in $B$ satisfying
\begin{itemize}
\item[(1)] for $a\in A$, $\rho(a)=0\Longrightarrow j_A(a)=0$,
\item[(2)] if $\rho(a)\in{({\psi}_V, \rho)}^{(1)}(\KK(M_L))$, then $j_A(a)\in j_A(K_{\alpha})$, 
\item[(3)] there exists a strongly continuous action $\beta:\T\to\Aut{\tau}_{\rho,V}(A{\rtimes}_{\alpha,L}\N)$ such that ${\beta}_z\circ{\tau}_{\rho,V}={\tau}_{\rho,V}\circ{\delta}_z$ for all $z\in\T$.
\end{itemize}
Then the corresponding representation ${\tau}_{\rho,V}:A{\rtimes}_{\alpha,L}\N\to B$ is faithful.
\end{thm}

The proof of Theorem~\ref{the giut} will use the following gauge-invariant uniqueness theorem for relative Cuntz-Pimsner algebras, which is due 
to Katsura \cite[Corollary~11.7]{k2} and Muhly-Tomforde \cite[\S5]{mt}. 

\begin{thm}\label{giut for rel CP algebras}
Suppose $X$ is a Hilbert bimodule over $A$ and $K$ is an ideal in $J(M_L)$. If $\mu:\OO(K,X)\to B$ is a homomorphism into a $C^*$-algbera $B$ satisfying
\begin{itemize}
\item[(i)] the restriction of $\mu$ to $k_A(A)$ is injective,
\item[(ii)] if $\mu(k_A(a))\in\mu\big({(k_X,k_A)}^{(1)}(\KK(X))\big)$, then $k_A(a)\in k_A(K)$,
\item[(iii)] there exists a strongly continuous action $\beta:\T\to\Aut\mu(\OO(K,X))$ such that ${\beta}_z\circ\mu=\mu\circ{{\gamma}}_z$ for all $z\in\T$,
\end{itemize} 
then $\mu$ is injective.
\end{thm}

\begin{proof}[Proof of Theorem~\ref{the giut}]
We will prove that ${\tau}_{\rho,V}\circ\theta$ satisfies the conditions of Theorem~\ref{giut for rel CP algebras}.
Suppose $a\in A$ satisfies $({\tau}_{\rho,V}\circ\theta)(k_A(a))=0$. ThenÊ 
\begin{align*}
\rho(a) = {\tau}_{\rho,V}(j_A(a))=Ê {\tau}_{\rho,V}(\theta(a))= 0,
\end{align*} 
which by (1) implies that $j_A(a)=0$. Hence $k_A(a)={\theta}^{-1}(j_A(a))=0$, and so ${\tau}_{\rho,V}\circ\theta$ is injective on $k_A(A)$.

Now suppose $a\in A$ and $({\tau}_{\rho,V}\circ\theta)(k_A(a))\in({\tau}_{\rho,V}\circ\theta)\big({(k_{M_L},k_A)}^{(1)}(\KK(M_L))\big)$. We have $({\tau}_{\rho,V}\circ\theta)(k_A(a))=\rho(a)$, and Lemma~\ref{equ for redundancies} gives
\begin{align*}
({\tau}_{\rho,V}\circ\theta)\big({(k_{M_L},k_A)}^{(1)}(\KK(M_L))\big) &= {\tau}_{\rho,V}\big({(\theta\circ k_{M_L},\theta\circ k_A)}^{(1)}(\KK(M_L))\big)\\
&= {\tau}_{\rho,V}\big({({\psi}_T,j_A)}^{(1)}(\KK(M_L))\big)\\
&= {\tau}_{\rho,V}\circ Q\big({({\psi}_S,i_A)}^{(1)}(\KK(M_L))\big)\\
&= (\rho\times V)\big({({\psi}_S,i_A)}^{(1)}(\KK(M_L))\big)\\
&= {({\psi}_V, \rho)}^{(1)}(\KK(M_L)).
\end{align*}
So $\rho(a)\in{({\psi}_V, \rho)}^{(1)}(\KK(M_L))$, and then it follows from (2) that $j_A(a)\in j_A(K_{\alpha})$. Hence $k_A(a)\in k_A(K_{\alpha})$.
By (3), we have
\[
\beta_z\circ \tau_{\rho,V}\circ \theta=\tau_{\rho,V}\circ\delta_z\circ \theta=\tau_{\rho,V}\circ \theta \circ\gamma_z,
\]
so Theorem~\ref{giut for rel CP algebras} implies that ${\tau}_{\rho,V}\circ\theta$ is injective. Thus ${\tau}_{\rho,V}$ is injective.
\end{proof} 

When the transfer operator $L$ is almost faithful on $K_{\alpha}$, our main theorem says that $j_A$ is injective. Using  \cite[Corollary~11.8]{k2} instead of Theorem~\ref{giut for rel CP algebras} yields the following gauge-invariant uniqueness theorem which directly generalises \cite[Theorem~4.2]{ev} (because the second condition ($2'$) trivially holds when $K_\alpha=J(M_L)$, as is the case when $\alpha(1)=1$).

\begin{cor}\label{the giut corollary}
Let $\alpha$ be an endomorphism of a unital $C^*$-algebra $A$, and let $L$ be a transfer operator for $(A,\alpha)$ which is almost faithful on $K_{\alpha}$. Suppose $B$ is a $C^*$-algebra and $(\rho,V)$ is a covariant representation of $(A,\alpha,L)$ in $B$ satisfying
\begin{itemize}
\item[($1'$)] $\rho$ is faithful,
\item[($2'$)] for $a\in J(M_L)$, $\rho(a)={({\psi}_V, \rho)}^{(1)}(\phi(a))$ implies $j_A(a)={({\psi}_T, j_A)}^{(1)}(\phi(a))$,
\item[(3)] there exists a strongly continuous action $\beta:\T\to\Aut{\tau}_{\rho,V}(A{\rtimes}_{\alpha,L}\N)$ such that ${\beta}_z\circ{\tau}_{\rho,V}={\tau}_{\rho,V}\circ{\gamma}_z$ for all $z\in\T$.
\end{itemize}
Then the corresponding representation ${\tau}_{\rho,V}:A{\rtimes}_{\alpha,L}\N\to B$ is faithful.
\end{cor}

\end{document}